\documentclass{article}[10pt]

\newcommand{\Div}{\mbox{\rm div}\,}

\newcommand{\grad}{\mbox{\rm grad}\,}
\newcommand{\supp}{\mbox{\rm supp}\,}

\newcommand{\Int}[2]{{\displaystyle \int_{ #1}^{ #2}}}

\newcommand{\Frac}[2]{\displaystyle{\frac{\displaystyle{#1}}{\displaystyle{#2}}}}

\newcommand{\beea}{\begin{eqnarray}}
\newcommand{\eeea}{\end{eqnarray}}
\newcommand{\ms}{\medskip\smallskip}

\newcommand{\BF}{\begin{footnotesize}}
\newcommand{\EF}{\end{footnotesize}}
\newcommand{\ode}[2]{{\displaystyle \frac{\mbox{$d #1$}}{\mbox{$d #2$}}}}

\newcommand{\bi}{\begin{itemize}}
\newcommand{\ei}{\end{itemize}}
\newcommand{\ed}{\end{document}}
\newcommand{\be}{\begin{equation}}
\newcommand{\ba}{\begin{array}}
\newcommand{\ea}{\end{array}}
\newcommand{\ee}{\end{equation}}
\newcommand{\eeq}[1]{\label{eq:#1}\end{equation}}
\newcommand{\real}{{\rm I\!\!\,R}}

\newcommand{\cala}{{\cal A}}
\newcommand{\calb}{{\cal B}}

\newcommand{\cald}{{\cal D}}

\newcommand{\calf}{{\cal F}}

\newcommand{\call}{{\cal L}}

\newcommand{\cals}{{\cal S}}

\newcommand{\calv}{{\cal V}}

\newcommand{\half}{\mbox{$\frac{1}{2}$}}

\def\Bbb R{\real}

\def\bar{\overline}

\newcommand{\ED}{\end{description}}

\newcommand{\Footnote}{~\footnote}

\newcommand{\Br}{\begin{remark}\begin{rm}}
\newcommand{\Er}{\end{rm}\end{remark}}
\newtheorem{remark}{Remark}[section]
\newcommand{\Bt}{\begin{theorem}\begin{sl}}
\newcommand{\Bd}{\begin{definition}\begin{sl}}
\newcommand{\Et}{\end{sl}\end{theorem}}
\newcommand{\Bl}{\begin{lemma}\begin{sl}}
\newcommand{\El}{\end{sl}\end{lemma}}
\newtheorem{theorem}{Theorem}[section]
\newtheorem{lemma}{Lemma}[section]
\newtheorem{definition}{Definition}[section]
\newtheorem{corollary}{Corollary}[section]
\newcommand{\Bc}{\begin{corollary}\begin{sl}}
\newcommand{\Ec}{\end{sl}\end{corollary}}
\newcommand{\EDD}[1]{\end{sl}\label{definition:#1}\end{definition}}
\newcommand{\ET}[1]{\end{sl}\label{theorem:#1}\end{theorem}}
\newcommand{\EL}[1]{\end{sl}\label{lemma:#1}\end{lemma}}
\newcommand{\theoref}[1]{{\rm Theorem \ref{theorem:#1}}}
\newcommand{\ER}[1]{\end{rm}\label{remark:#1}\end{remark}}
\newcommand{\EC}[1]{\end{sl}\label{corollary:#1}\end{corollary}}

\newcommand{\lemmref}[1]{{\rm Lemma \ref{lemma:#1}}}
\usepackage{amssymb}
\usepackage{amsmath}
\usepackage[mathscr]{eucal}
\newcommand{\ds}{\displaystyle}

\newcommand{\plq}{\rule{0.in}{.25in}\left[}
\newcommand{\prq}{\rule{0.in}{.25in}\right]}

\renewcommand{\grad}{{\nabla}}

\newcommand{\ints}{\int_{\Sigma}}
\newcommand{\gd}{\grad}
\newcommand{\dv}{{\Div }}

\newcommand{\cd}{{\cal{D}}}

\newcommand{\mr}{{\real}}

\newcommand{\id}{{\textrm{I}}}
\renewcommand{\real}{{\mathbb R}}
\begin{document}

\title{Large-Time Behavior of a Rigid Body of Arbitrary Shape in a Viscous Fluid Under the Action of Prescribed Forces and Torques} 
        
\author{Giovanni P. Galdi\Footnote{Department of Mechanical Engineering and Materials Science, University of Pittsburgh, Pittsburgh, USA.}}
\date{\it\small In loving memory of Olga~Ladyzhenskaya, a Founder of mathematical fluid mechanics.}          
\maketitle
\begin{abstract} 
Let $\calb$ be a sufficiently smooth rigid body (compact set of $\real^3$) of arbitrary shape moving in an unbounded Navier-Stokes liquid  under the action of prescribed external force, $\textup{F}$, and torque, $\textup{M}$. We show that if the data are suitably regular and small, and $\textup{F}$ and $\textup{M}$ vanish for large times in the $L^2$-sense, there exists at least one global strong solution to the corresponding initial-boundary value problem. Moreover, this solution converges to zero  as time approaches infinity. This type of results was known, so far, only when $\calb$ is a ball.  
\end{abstract}
\renewcommand{\theequation}{\arabic{section}.\arabic{equation}}\setcounter{equation}{0}
\section{Introduction}  The motion of a (finite) rigid body, $\calb$, in an unbounded Navier-Stokes liquid has been the object of a number of deep researches. Particularly intriguing is the case when the motion of the body is not given and, in general,  one prescribes total force, $\textup{F}$, and torque $\textup{M}$, acting on it. Since the 
presence of the body {affects the flow of the liquid},  and this, in 
turn, { affects the motion of the body}, the problem of determining the 
flow characteristics thus becomes  highly coupled. It is  this distinctive  property that makes  any  mathematical problem related to body-liquid interaction especially interesting and challenging. 
\par    
In this paper we are interested in the study of two basic questions related to the situation just described when the shape of $\calb$ is {\em not} specified, and precisely: (i) existence of a global-in-time strong solutions to the relevant initial-boundary value problem, and (ii) their asymptotic behavior for all large times. Before stating our results, we would like to recall all known major contributions related to this type of investigation, which  will also furnish the motivation for the present study.\footnote{We shall restrict ourselves to the three-dimensional case that is the focus of our work.} The first existence result is due to Serre \cite{ser87}, who proves global existence of weak solutions a la Leray-Hopf. As the author himself observes, the proof is exactly the same as the classical one for the  Navier-Stokes problem and presents no challenges. Instead, a less obvious task is  to show existence of  {\em strong} solutions  having enough regularity  as to solve the given equations (at least) at almost every point of the space-time. This question was first successfully tackled by Galdi \& Silvestre \cite{GaSi} who proved existence of strong  solutions, in the sense of Prodi-Ladyzhenskaya \cite{Lad,Pro}, for data of arbitrary ``size" --in a suitable class-- at least in a time-interval $[0,T)$ for some $T>0$. Successively, a similar result, but with a different approach, was established by Cumsille \& Tucsnak \cite{CuTu}, when the body is  allowed to rotate but not to translate, and $\textup{F}\equiv\textup{M}\equiv 0$.  In both papers \cite{CuTu,GaSi} the functional framework is the $L^2$ Hilbert-setting. The study of existence of strong solutions in the $L^q$ setting, $q\in (1,\infty)$, was initiated by Wang \& Xin \cite{WaXi}, who established a local in time result with $\textup{F}\equiv\textup{M}\equiv 0$, in the special case when $\calb$ is a ball. 
Since the general shape of the body is a most relevant feature of our result here, let us briefly  comment on how the hypothesis of $\calb$ being a  ball brings in some basic simplifications and mathematical  properties that are lost in the general case. In the first place, this assumption  eliminates the presence, in the linear momentum equation, of a term whose coefficient becomes unbounded at large spatial distances. Furthermore, as shown in \cite{CuTa,WaXi}, the relevant linear operator, suitably defined, is the generator of an analytic semigroup, a property that is no longer valid for bodies of arbitrary shape, just because of the occurrence of the unbounded term \cite{HiSa,Sh}.      
Local $L^q$ existence for $\calb$ of {\em arbitrary} shape was successively established by Geissert {\em et al.} \cite{GGH}, by  maximal regularity theory, again with $\textup{F}\equiv\textup{M}\equiv 0$. Concerning the question of {\em global} existence, it was first studied and positively answered by Cumsille \& Takahashi \cite{CuTa}. In particular, they showed that if, in appropriate norms, the initial data are ``small"  and the external forces are summable over the whole half-line $(0,\infty)$ and ``small,"   there exists a (unique) corresponding solution defined for all times and belonging to a functional class similar to that considered in \cite{GaSi}. The method used in \cite{CuTa} is based on a particular cut-off technique that, on one hand, eliminates the difficulty due to the unbounded coefficient, but, on the other hand, is not able to provide any information on the large-time behavior of the solutions that, under the given assumption, are  expected to reach, eventually, the rest-state. The question of the asymptotic behavior of solutions (along with their global existence) has been analyzed very recently  in a remarkable paper by Ervedoza {\em   
et al.} \cite{EMT}, when  $\textup{F}\equiv\textup{M}\equiv 0$. The main tool
is new $L^p- L^q$ estimates for the 
body-liquid semigroup. Even though the estimates are proved
for bodies of arbitrary shape, their use in showing global existence of solutions (for small data) along with their asymptotic decay to  rest-state requires $\calb$ to be a ball. In such a case, the authors also provide a sharp decay rate  that implies that the center of mass of $\calb$ can only cover a finite distance from its initial position, as expected on physical grounds.  
\par
In view of all the above, the following basic question --brought to my attention by Professor Toshiaki Hishida--  remains still open: Let $\calb$ be of {\em arbitrary} shape, subject to prescribed force and torque that vanish (in suitable sense) as time goes to infinity. Does the body-liquid problem have a global solution that, in addition, ultimately tends to the rest-state?  
\par
Objective of this note is to give a positive answer to this question, under the assumption of ``small" data. More precisely, we  show (see \theoref{nuovo}) that the local solution constructed  in \cite{GaSi} can, in fact, be extended to arbitrarily positive times, if the data are small enough. This is accomplished by means of a generalized Gronwall's lemma, proved in \lemmref{4.3}. By the same tool and under the same hypotheses, we then prove that solutions must eventually converge to the state of rest. Unfortunately, we are not able to furnish a rate of decay, which thus leaves room to further investigation. However, in the case when  $\calb$ is a ball, we do provide such a decay that, in the $L^2$ framework considered here, appears to be rather sharp; see Remark \ref{rem:4.1}. 
\par
The plan of the paper is as follows. In Section 2 we formulate  the problem and state our main result in  \theoref{nuovo}. Successively,  in  Section 3 we prove two  basic ``energy equations" valid in the class of solutions considered in \cite{GaSi}; see \lemmref{3.3}. This requires some estimates on the time derivative of the velocity field and on the pressure field that are carried out in \lemmref{3.1} and \lemmref{3.2}. Finally, in Section 4 we give a proof of \theoref{nuovo}, by combining estimates obtained from the energy equations with the Gronwall-like lemma showed in \lemmref{4.3}.

\section{Mathematical Formulation and Main Result}
A rigid body ${\cal B}$ --that is, a sufficiently smooth, compact and connected set of $\real^3$-- is fully immersed in  a quiescent Navier-Stokes liquid, ${\cal L}$,  that fills the entire three-dimensional space exterior to $\calb$. We suppose that, with respect to an inertial frame, $\mathscr F$,  the body  is subject to prescribed force, $\textup{F}=\textup{F}(t)$, and torque, $\textup{M}=\textup{M}(t)$, $t\ge0$.  Following a standard procedure, we shall describe the motion of the coupled system $\cals:=\{\calb,\call\}$ with respect to a frame, $\mathscr S$, attached to $\calb$ and with its origin at the center of mass, $G$, of $\calb$. In such a way, in particular, the domain occupied by $\call$ becomes time-independent, and we will  denote it by $\cald$ ($:=\real^3\backslash\calb$) and  by $\Sigma$ its boundary. 
Assuming, without loss of generality, that $\mathscr S$ and $\mathscr F$ coincide at $t=0$,  the equations governing the motion of $\cals$  in $\mathscr S$ are given by (see \cite{GaRev})
\begin{equation}
\left.
\begin{array}{c}\ms
\left. 
\ds
\begin{array}{c}\ms
\varrho \partial_t{u} = {\Div}T(u,p) - \varrho  [ (u-V)\cdot \grad u 
 + \omega \times u ]
\\ 
{\Div }u=0
\end{array}
\ds
\right\} {\rm in }\,\ \cd \times (0,\infty ) \\ \ms
\ds 
u = V  \ \ {\rm at }\ \Sigma \times (0,\infty ) \\ \ms
\ds
\lim_{|x| \rightarrow \infty }u(x,t)=0,{ }\,\ t\in
(0,\infty ) \\ \ms 
\ds 
m\dot{\xi }+ m \omega \times \xi + \int_{\Sigma } T(u,p)\cdot n = \textsf{F} \\ 
\ms \ds
{\sf I}\cdot \dot{\omega }+  \omega \times ({\sf I} \cdot \omega )+ \int_{\Sigma }x \times T(u,p)\cdot n =  \textsf{M} \\
\ms
\ds 
\xi (0)=\xi _{0},\quad \ \omega (0)=\omega _{0} \\
\ds
u(x,0) = u_{0}(x),{ }\, x\in \mathcal{D}\,.
\end{array}
\label{eqpmodd}
\right.
\end{equation} 
Here, $u$ and $p$ are  velocity and pressure fields of ${\cal L}$, $\varrho$ its (constant) density, and $V(x,t):= \xi(t) + \omega (t) \times x,$ where $\xi$ is the velocity of the center of mass of $\calb$ and $\omega$ its angular velocity. Also, $T$ is the Cauchy stress tensor given by
$$
T(u,p)=2\mu\,D(u)-p\,\id,\ \ D(u):=\half\, \big(\grad u+(\grad u)^\top\big)\,,
$$
with $\mu$ shear-viscosity coefficient and $\id$ identity.
Moreover,  $m$  is the mass of  $\calb$ and ${\sf I}$ its inertia tensor relative to $G$. Furthermore,
\begin{equation}
\left\{
\begin{array}{c}\ms
\textsf{F}(t)=Q^{\top }(t)\cdot\textup{F}(t), \\
\textsf{M}(t)=Q^{\top }(t)\cdot\textup{M}(t)\,,
\end{array}
\right.
\label{eqfm}
\end{equation}
with the tensor $Q$ satisfying the following equation
\be\label{Q}
\left\{ 
\begin{array}{c}\ms
\ds \dot{Q} = -Q\cdot \Omega(\omega)\\ \ds Q(0)=\id  
\end{array}
\right.
\qquad \qquad \qquad \Omega(\omega)= \plq \begin{array}{ccc}0&\omega_3&-\omega_2\\-\omega_3&0&\omega_1\\\omega_2&-\omega_1&0 \end{array}\prq
\end{equation}
In particular, $Q$ is proper orthogonal, that is,
$$
\ds Q^{\top}(t)\cdot Q(t)=Q(t)\cdot Q^{\top}(t)=\id, \ \ \
\ds \det Q(t)=1, \quad \mbox{for all}\, t \in \real\,.
$$
\par
In order to state our main result, we need a suitable function space.
Let 
$$
{\cal R}:=\{\overline{u}\in C^{\infty}(\mr^3): \overline{u}(x)=\overline{u}_1 + \overline{u}_2 \times x, \quad \overline{u}_1,\overline{u}_2 \in \mr^3\}\,,
$$
and define\Footnote{
We shall use standard notation for function
spaces, see \cite{ad}. So, for instance, 
$L^q(\cala)$, $W^{m,q}(\cala)$, $W_0^{m,q}(\cala)$,
etc., will denote the usual Lebesgue and Sobolev spaces on the
domain $\cala$, with norms
$\|\,\cdot\|_{q,\cala}$ and $\|\,\cdot\|_{m,q,\cala}$, respectively.
Whenever confusion will not arise, we shall omit the subscript $\cala$.
The trace space on $\partial\cala$ for functions from
$W^{m,q}(\cala)$ will be denoted by $W^{m-1/q,q}(\partial\cala)$ and its norm
by $\|\,\cdot\|_{m-1/q,q,\partial\cala}$.
Occasionally, for $X$ a Banach space, we denote by $\|\cdot\|_X$ its associated norm. Moreover $L^q(I;X)$, $C(I;X)$ $I$  real interval,  denote classical Bochner spaces.} 
$$
{\cal V}(\cd) =   \{ u \in W^{1,2}(\cd): \dv u =0 \textrm{ in }\cd, \textrm{ }\textrm{ } u\left|_{\Sigma}\right.= \overline{u}, \ \mbox{for some} \  \overline{u}   \in {\cal R}\}.
$$
We also set
$$\ba{ll}\medskip
B_R:=\{x\in\real^3:\,|x|< R\}\,;\ \ R_*:= 2\inf\,\{R\in (0,\infty): \calb\cap B_R\supset\calb\}\,;\\ \cald_R:=\cald\cap B_R\,,\ \ \cald^R=\cald\backslash \bar{\cald_R}\,,\ \ R>R_*\,.\ea
$$
\par
The main objective of this paper is to show the following result.
\Bt
Let ${\cald}$ be of class $C^2.$ Let $\textup{F},\textup{M}\in L^2(0,\infty)$ and  $u_0 \in {\cal V}(\cd)$ with $u_0|_{\Sigma}=\xi_0+\omega_0\times x$. Then, there is $\delta>0$ such that if
\be\label{smalldata}
\|u_0\|_{1,2}+|\xi_0|+|\omega_0|+\|\textup{F}\|_{L^2(0,\infty)}+\|\textup{M}\|_{L^2(0,\infty)}\le\delta\,,
\ee
there exist functions $u=u(x,t),$ $p=p(x,t),$ $\xi=\xi(t),$ $\omega=\omega(t)$, and $Q=Q(t)$ satisfying {\rm (\ref{eqpmodd})-(\ref{Q})} {a.e.}, such that 
\be
\begin{array}{ll}\medskip
u \in L^{\infty}(0,\infty;W^{1,2}(\cd)),\ \ \
\gd u \in L^2(0,\infty;W^{1,2}(\cd))\\\medskip
\xi,\textrm{ } \omega\in W^{1,2}(0,\infty), \ \ \ \nabla p\in L^2(0,\infty;L^2(\cald))\,, \ \ \ Q\in W^{2,2}(0,\infty)\\
\partial_t u,\textrm{ } p\in L^2(0,\infty;L^2(\cd_R)), \ \mbox{for all}\, R \ge R_*\,.
\end{array}
\eeq{class}
Moreover, for all $T>0$,
\begin{eqnarray*}
&&\xi,\textrm{ } \omega,\, Q \in C([0,T]), \ \textrm{ with }\xi(0)=\xi_0,\textrm{ } \omega(0)=\omega_0,\textrm{ } Q(0)={\rm I}\\
&&u \in C([0,T];W^{1,2}(\cd_R)), \ \mbox{for all}\, R \ge R_*, \ \textrm{ with } u(.,0)=u_0(.).
\end{eqnarray*}
Finally,
\be
\lim_{t\to\infty}\big(\|u(t)\|_6+\|\nabla u(t)\|_2+|\xi(t)|+|\omega(t)|\big)=0\,.\label{asym}
\ee
\ET{nuovo}
\par
Before carrying out, in the following sections, the proof of the theorem,  we would like to make some comments. The major aspect of our results is expressed by the asymptotic property (\ref{asym}), which states that, eventually, the coupled system $\cals$ will go to rest, independently of the shape of $\calb$. In fact, to date, this property was known only when $\calb$ is a ball \cite{EMT} and $\textup{F}\equiv\textup{M}\equiv 0$. However, unlike \cite{EMT}, even with the additional assumption on $\textup{F}$ and $\textup{M}$, we are not able to furnish a rate of decay. We may guess that it is $O(t^{-\frac12})$, but a proof seems to be currently out of reach; see also Remark \ref{rem:4.1}. \par On the other side, in the case when either $\textup{F}$ or $\textup{M}$ is time-independent the existence of global  strong solutions and, more intriguingly, the assessment of their asymptotic behavior represents   a formidable open question. A remarkable example is the free-falling body problem where $\textup{M}=0$ and $\textup{F}=m_e g$, with $m_e$ buoyant mass of $\calb$ and $g$ acceleration of gravity. In such a case it is expected that, at least for small $m_e$, the coupled system $\cals$ will tend, as $t\to\infty$, to a steady-state configuration. However,  as shown in \cite{GaVa}, the steady-state problem may have multiple solutions, even for vanishingly small $m_e$. One may thus conjecture that $\cals$ will approach, eventually,  one of the  locally unique, stable configuration that are experimentally observed, at least when $\calb$ has fore-and-aft symmetry, like homogeneous cylinder \cite{JL}. Nevertheless, even in this case,  a rigorous proof is far from obvious.\footnote{For a formal proof, see \cite{Cox}.}   

\section{Preliminary Results}
The goal of this section is to derive a number of {\em a priori} estimates for solutions to \eqref{eqpmodd}--\eqref{Q} in a suitable function class that we define next.
\Bd We say that $(u,p,\xi,\omega,Q)$ is in the class $\mathscr C_T$,   some $T\in(0,\infty]$, if, for all $\tau\in (0,T)$,
$$
\begin{array}{ll}\medskip
u \in L^{\infty}(0,\tau;W^{1,2}(\cd)),\ \ \
\gd u \in L^2(0,\tau;W^{1,2}(\cd))\\\medskip
\xi,\textrm{ } \omega\in W^{1,2}(0,\tau), \ \ \ \nabla p\in L^2(0,\tau;L^2(\cald))\,, \ \ \ Q\in W^{2,2}(0,\tau)\\
u \in C([0,\tau];W^{1,2}(\cd_R)) \,,\ \ \partial_t u,\textrm{ } p\in L^2(0,\tau;L^2(\cd_R))\,,\ \ \mbox{for all}\, R \ge R_*\,.
\end{array}
$$
\EDD{d}
The following results hold.
\Bl Let $(u,p,\xi,\omega,Q)$ be  a solution to  \eqref{eqpmodd}--\eqref{Q} in the class $\mathscr C_T$. Then for a.a. $t\in(0,T)$
$$
\partial_tu/r\in L^2(\cald)\,,\ \ r:=(x_ix_i)^{\frac12}\,.
$$
\EL{3.1}
{\em Proof.} From the assumption, we immediately show
\be\label{cat1}
\big(\Div T(u,p)+\varrho\,(V\cdot\nabla u-\omega\times u\big)/r\in L^2(\cald)\,.
\ee
Moreover, by Schwarz and Sobolev inequalities
\be\label{cat2}
\|u\cdot\nabla u\|_2\le \|u\|_4\|\nabla u\|_4\le c\,\|u\|_{2,2}^2\,.
\ee
Thus, since $r^{-1}\in L^\infty(\cald)$, the lemma follows from \eqref{cat1}, \eqref{cat2} and  \eqref{eqpmodd}$_1$.
\par\hfill$\square$\par
\Bl Let $(u,p,\xi,\omega,Q)$ be  a solution to  \eqref{eqpmodd}--\eqref{Q} in the class $\mathscr C_T$. Then for a.a. $t\in(0,T)$
$$
\nabla p\in L^{q_1}(\cald^{2R_*})\,,\ \ p\in L^{q_2}(\cald^{2R_*})\,,\ \ \mbox{for all $q_1\in (1,6]$\,,\ $q_2\in (\frac32,\infty]$}\,.
$$
\EL{3.2}
{\em Proof.} Observing that, in the sense of distribution, 
$$
\Div[\varrho\big(\partial u_t-(V-u)\cdot\nabla u\big)-\mu\Delta u]=0\,,
$$
from \eqref{eqpmodd}$_1$ we get for a.a. $t\in (0,T)$
\be\label{cat3}
\Delta p=\Div f\,\ \textrm{in}\ \cald\,,\ \ f:=\varrho\,u\cdot\nabla u\,.
\ee
Let $\psi=\psi(|x|)$ be a smooth, non-decreasing function such that $\psi(|x|)=0$, if $|x|\le R_*)$, while $\psi(|x|)=1$ if $|x|\ge 2R_*$. 
Setting ${\sf p}:=\psi\, p$, and extending ${\sf p}$ to zero outside $\cald^{R_*}$, from \eqref{cat3} we find
\be
\Delta {\sf p}=F\ \ \mbox{in $\real^3$}\,,\label{Poi} 
\ee
where
\be
F:=\Div(2 p\,\nabla\psi+\psi\,f)-\nabla\psi \cdot f-p\,\Delta\psi\,.
\label{sfax}
\ee
Since $u,p\in\mathscr C_T$, from Sobolev embedding theorem it follows that
$$
u\in L^q(\cald)\,,\ \ \nabla u\in L^s(\cald)\,,\ \ p\in L^r_{\rm loc}(\bar{\cald})\,,\ \ \mbox{for all $q\in [2,\infty)$, $s\in [2,6]$, $r\in [1,6]$\,.}
$$
This implies, in particular,
$$
f\in L^r(\cald)\,,\ \ \mbox{for all $r\in [1,6]$}\,.
$$
Therefore,
$$
\|F\|_{-1,r}\le c\,\left(\|f\|_{r}+\|p\|_{r,\cald_{2R_*}}\right)\,,\ \ \mbox{for all $r\in [1,6]$}
$$ 
Problem \eqref{Poi}--\eqref{sfax}  formally coincides with problem (III.1.20) studied in \cite[pp. 149-150]{GaB}, for which, observing that  $\nabla {\sf p}\in L^2(\real^3)$ and recalling that  $\psi_R\equiv 1$ in $\cald^{2R_*}$,  one proves (see \cite[eq. (III.1.23)]{GaB}) that $\nabla p\in L^{r}(\cald^{2R_*})$ for all $r\in (1,6]$. By \cite[Theorems II.2.1(i) and II.9.1]{GaB} the latter in turn implies $p\in L^\sigma(\cald^{{2R_*}})$, for all $\sigma\in (\frac32,\infty]$, which completes the proof.\par\hfill$\square$\par
\Bl Let $(u,p,\xi,\omega,Q)$ be  a solution to  \eqref{eqpmodd}--\eqref{Q} in the class $\mathscr C_T$. Then, the following relations hold, for all $t\in(0,T)$
\be\ba{ll}\ms 
\frac 12\Frac d{dt}\bigg( \varrho  \| u \|_2^2 + m \left| \xi
\right|^2 + \omega \cdot I \cdot \omega  \bigg) + 2 \mu\| D(u) \|_2^2={\sf F}\cdot\xi+ {\sf M}\cdot \omega\\\ms
\ds \mu\frac{d}{dt}\| D(u)\|_{2}^2  +  m |\dot{\xi}|^2 + \dot{\omega} \cdot {\sf I} \cdot \dot{\omega} + \| \dv T(u,p)\|_{2}^2 
= -m \omega \times \xi \cdot \dot{\xi}  -  \omega \times ({\sf I} \cdot \omega )\cdot \dot{\omega}- {\sf F} \cdot \dot{\xi} - {\sf M} \cdot \dot{\omega} \\
\quad\,+    \varrho\ds\int_{\cd} u\cdot \grad u \cdot \dv T(u,p)
+  \varrho\,\mu\ds\int_{\cd} ( \omega \times \grad u_i\cdot\nabla u_i - \nabla(\omega \times u):\nabla u )
-\varrho\,\mu\Int{\Sigma}{}\left(n\cdot\nabla u\cdot\Phi-\half V\cdot n|\nabla u|^2\right),
\end{array}\label{casa}
\end{equation}
where $\Phi:=V\cdot\nabla u-\omega\times u$\,.
\EL{3.3}
{\em Proof.} Let 
$$
\ds T_1:= \frac{1}{m} \ints  T(u,p)\cdot n, \ \ T_2:= {\sf I}^{-1}\cdot \ints x \times T(u,p)\cdot n.
$$ 
We test both sides of \eqref{eqpmodd}$_1$ by $u$,  integrate by parts over $\cald_R:=\cald\cap \{|x|<R\}$, $R>R_*$, and use \eqref{eqpmodd}$_{2,3}$ to get
\be
\half\varrho\ode{}t\|u(t)\|_{2,\cald_R}^2+2\mu\|D(u)\|^2_{2,\cald_R}=m\xi\cdot T_1+\omega\cdot{\sf I}\cdot T_2 +\int_{\partial B_R}\left(u\cdot T(u,p)\cdot n-\half \varrho u^2(u-\xi)\cdot n\right)\,,
\label{brodo1}
\ee
where we observed that $\omega\times x\cdot n=0$ at $\partial B_R$. Since $u\in \mathscr C_T$, and also with the help of \lemmref{3.2}, it is readily seen that the surface integral in \eqref{brodo1}  is in $L^1(R_*,\infty)$, so that we may let $R\to\infty$ along a sequence to get 
\be
\half\varrho\ode{}t\|u(t)\|_{2}^2+2\mu\|D(u)\|^2_{2}=m\xi\cdot T_1+\omega\cdot{\sf I}\cdot T_2 \,.
\label{brodo2}
\ee
As a result, employing \eqref{eqpmodd}$_{5,6}$
 in \eqref{brodo2} we deduce \eqref{casa}$_1$. In order to show  \eqref{casa}$_2$, we begin to observe that \cite[Lemma 2.4(b)]{GaSi}
\be
\label{phin}
\Phi\cdot n=0\ \ \mbox{at $\Sigma$}\,.
\ee
Moreover, for any $R>R_*$, let $\psi_R=\psi_R(|x|)$ be a non-decreasing, smooth function such that $\psi_R(|x|)=1$, if $|x|\le R$ and $\psi_R(|x|)=0$, if $|x|\ge 2R$, and
\be\label{psi}
|\nabla \psi_R(|x|)|\le C\,R^{-1}\,,
\ee
with $C$ independent of $x$ and $R$.
We next test  \eqref{eqpmodd}$_1$ by $\psi_R\,\Div T(u,p)$ to get
\be\label{KTM}
\int_\cald\psi_R\partial u_t\cdot\Div T=\|\sqrt{\psi_R}\Div T\|_2^2-\varrho\int_\cald\left(\psi_Ru\cdot\nabla u\cdot\Div T-\psi_R\Phi\cdot\Div T\right)\,. 
\ee
By integration by parts, we show
$$\ba{rl}\medskip\ds
\int_\cald\psi_R\partial u_t\cdot\Div T&=\ds \int_\cald\left[\Div(\psi_R\,\partial_tu\cdot T)-2\mu\,\psi_R\,D(\partial_tu):D(u)\right]\\&=\ds \int_\Sigma \dot{V}\cdot T\cdot n-\mu\ode{}t \|\sqrt{\psi_R}D(u)\|_2^2-\int_\cald \nabla\psi_R\cdot T\cdot\partial_tu\,. 
\ea
$$
Using  \eqref{eqpmodd}$_{5,6}$ in the surface integral, we deduce
\be\label{KTM1}
\ba{rl}\medskip\ds
\int_\cald\psi_R\partial u_t\cdot\Div T=&\!\!-m\dot{\xi}^2-\dot{\omega}\cdot{\sf I}\cdot\dot{\omega}-m\omega\times\xi\cdot\dot{\xi}-\omega\times({\sf I}\cdot\omega)\dot{\omega}-{\sf F}\cdot\dot{\xi}-{\sf M}\cdot\dot{\omega}\\&\ds -\mu\ode{}t \|\sqrt{\psi_R}D(u)\|_2^2-\int_\cald \nabla\psi_R\cdot T\cdot\partial_tu\,. 
\ea
\ee
Next, integrating by parts and with the help of \eqref{phin} we show
\be\label{KTM2}
\Int{\cald}{}\psi_R\Phi\cdot\Div T=2\mu\Int{\Sigma}{}\Phi\cdot D(u)\cdot n-2\mu\Int{\cald}{}\psi_R\partial_i\Phi_jD_{ij}(u)-\Int{\cald}{}\nabla\psi_R\cdot T\cdot\Phi\,.
\ee
Now, using $\Div u=\Div V=0$,
$$\ba{rl}\medskip
2\partial_i\Phi_jD_{ij}(u)&=\partial_i(\Phi_j\partial_j u_i)+(\partial_iV)\cdot\nabla u_j\partial_iu_j+\half V\cdot\nabla(|\nabla u|^2)-\nabla(\omega\times u):\nabla u\\
&=\Div(\Phi\cdot\nabla u+\half V\,|\nabla u|^2)+\omega\times \nabla u_i\cdot\nabla u_i-\nabla(\omega\times u):\nabla u\,.
\ea 
$$
Substituting the latter in \eqref{KTM2} and using Gauss theorem, we infer
\be\ba{ll}\medskip\label{KTM3}
\Int{\cald}{}\psi_R\Phi\cdot\Div T=&\!\!
-\mu\Int{\cald}{}\psi_R\left(\omega\times \nabla u_i\cdot\nabla u_i-\nabla(\omega\times u):\nabla u\right)+\mu\Int{\Sigma}{}\left(n\cdot\grad u\cdot\Phi-\half V\cdot n|\nabla u|^2\right)\\
&+\Int{\cald}{}\nabla\psi_R\cdot\left[2\mu(\nabla u^\top\cdot\Phi+\half V|\nabla u|^2)-T\cdot\Phi\right]\,.
\ea
\ee
Collecting \eqref{KTM}, \eqref{KTM1} and \eqref{KTM3} we deduce
\be\label{KTM4}\ba{ll}\medskip
\mu\ode{}t \|\sqrt{\psi_R}D(u)\|_2^2
+\|\sqrt{\psi_R}\Div T\|_2^2+m\dot{\xi}^2+\dot{\omega}\cdot{\sf I}\cdot\dot{\omega}=-m\omega\times\xi\cdot\dot{\xi}-\omega\times({\sf I}\cdot\omega)\,\dot{\omega}-{\sf F}\cdot\dot{\xi}-{\sf M}\cdot\dot{\omega}\\\medskip \ds +\varrho\Int{\cald}{}\psi_R\left[u\cdot\nabla u\cdot\Div T +\mu\left(\omega\times \nabla u_i\cdot\nabla u_i-\nabla(\omega\times u):\nabla u\right)\right]-\varrho\mu\Int{\Sigma}{}\left(n\cdot\grad u\cdot\Phi-\half V\cdot n|\nabla u|^2\right)\\
-\varrho\Int{\cald}{}\nabla\psi_R\cdot\left[T\cdot\partial_tu+2\mu(\nabla u^\top\cdot\Phi+\half V|\nabla u|^2)-T\cdot\Phi\right]
\ea
\ee
Let us denote by $\calf(R)$ the last integral on the right-hand side of \eqref{KTM4}. Observing that $\supp(\psi_R)\subseteq \{R\le |x|\le 2R\}$ and recalling \eqref{psi}, we obtain
$$
|\calf(R)|\le c\Int{R\le |x|\le 2R}{}|x|^{-1}\left[(|T|(|\partial_tu|+|\Phi|)+|\nabla u||\Phi|+|V||\nabla u|^2\right]\le c \Int{\cald^R}{}\left(|x|^{-2}|\partial_tu|^2+|\nabla u|^2+|u|^2+|p|^2\right).
$$
As a result, by assumption,  \lemmref{3.1}, and \lemmref{3.2}, we show
\be\label{KTM5}
\lim_{R\to\infty}{\mathcal F}(R)=0\,.
\ee
Furthermore, since $(u,p)\in\mathscr C_T$ it readily checked that
\be\label{KTM6}
\left[u\cdot\nabla u\cdot\Div T +\mu\left(\omega\times \nabla u_i\cdot\nabla u_i-\nabla(\omega\times u):\nabla u\right)\right]\in L^1(\cald)\,.
\ee
We integrate both sides of \eqref{KTM3} over $(0,t)$, $t\in (0,\tau]$, let $R\to\infty$ and employ \eqref{KTM5}, \eqref{KTM6} along with Lebesgue dominated convergence theorem. If we differentiate with respect to $t$ the resulting equation, we then end up with \eqref{casa}$_2$, which completes the proof of the lemma.\par\hfill$\square$\par  \setcounter{equation}{0} 
\section{Proof of the Main Result}
In this section we shall prove  \theoref{nuovo}. Before beginning the proof, however, we premise  further lemmas. 
\Bl Let $u\in \calv(\cal D)$. Then, 
$$
\|\nabla u\|_2=\sqrt{2}\|D(u)\|_2\,.
$$
Moreover, 
there is $c_1=c_1(\cald)$ such that
$$
|\bar{u}_1|+|\bar{u}_2|\le c\,\|D(u)\|_2\,.
$$
Finally, $u\in L^6(\cald)$ and there is $c_2=c_2(\cald)$ such that
$$
\|u\|_6\le c_2\,\|D(u)\|_2\,.
$$
\EL{4.1}
{\em Proof.} See \cite[Section 4.2.1]{GaRev}.\hfill$\square$
\Bl Let $u\in \calv({\cal D})\cap W^{2,2}(\cald)$, $\nabla p\in L^2(\cald)$. Then, there is $C=C(\cald)$ such that
$$
\|D^2u\|_2\le C\,(\|\Div T\|_2+\|D(u)\|_2)\,.
$$
\EL{4.2}
{\em Proof.} The lemma follows from \cite[Lemma V.4.3]{GaB} and \lemmref{4.1}.\hfill$\square$
\Bl Let $y:(0,T)\mapsto [0,\infty)$, $T>0$, be absolutely continuous,  such that 
\be
y'(t)\le G(t)+c_1y(t)+c_2y^\alpha(t)\,,\ \ \alpha>1,\ \ t\in(0,T)\,,
\label{1}
\ee
where $G\in L^1(0,T)$, $G(t)\ge 0$, and $c_i\in [0,\infty)$, $i=1,2$.
Then, if --in case at least one of the constants $c_i$ is not zero-- it is also $y\in L^1(0,T)$, there exists $\eta>0$, such that from
\be
y(0)+\int_0^T G(s)\,ds+\int_0^T y(s)\,ds\le\eta\label{2}
\ee
it follows $y\in L^\infty(0,T)$ and
\be\label{Co1}
y(t)< M\,\eta\,,\ \   M:=2\max\{1,c_1,c_2\}, \ \ t\in [0,T)\,.
\ee
In the case $T=\infty$, we also have
\be\label{Co2}
\lim_{t\to\infty}y(t)=0\,, 
\ee
and, if $t\,G\in L^1(0,\infty)$, $c_1=0$, and $\alpha\ge2$, even
\be\label{Co3}
\sup_{t\in (0,\infty)}\left(t\,y(t)\right)\le A<\infty\,.
\ee
\EL{4.3}
{\em Proof.} Since $y(0)\le\eta $, contradicting \eqref{Co1} means that there  exists $t_0\in (0,T)$ such that $y(t)<M\,\eta$ for all $t\in [0,t_0)$ and $y(t_0)=M\,\eta$. Integrating both sides of \eqref{1} from 0 to $t_0$, we deduce, in particular
$$
y(t_0)\le y(0) +\int_0^{T}G(s)\,{\rm d}s+c_1\int_0^T y(s)\,{\rm d}s +c_2\int_0^{t_0}y^\alpha(s)\,{\rm d}s\,.  
$$
Therefore, setting $\mu=\max\{1,c_1,c_2\}$, from this inequality and \eqref{2} we find
$$
y(t_0)\le \mu\,\eta\left[1 + (M\eta)^{\alpha-1}\right]
\,,
$$ 
so that, choosing $\eta\in (0,(1/M)^{1/(\alpha-1)})$, we obtain $y(t_0)<M\,\eta$, a contradiction that proves \eqref{Co1}.
In order to show the property \eqref{Co2}, we observe that, in view of assumption \eqref{2} and being $T=\infty$, there exists an unbounded  sequence $\{t_k\}\subset (0,\infty)$ such that
\be\label{3}
\lim_{k\to\infty}y(t_k)=0
\,.
\ee
We then integrate \eqref{1} from $t_k$ to arbitrary $t>t_k$,  and recall \eqref{Co1}, to deduce, in particular,
$$
y(t)\le y(t_k) +\int_{t_k}^{\infty}G(s)\,{\rm d}s+c_3\int_{t_k}^\infty y(s)\,{\rm d}s\,,  
$$
for some $c_3>0$ and all $t>t_k$. In view of \eqref{2} with $T=\infty$ and  \eqref{3}, the right-hand side of this inequality can be made as small as we please, by taking sufficiently large $k$, and property \eqref{Co2} follows. Finally, take in \eqref{1} $c_1=0$ and $\alpha\ge2$. Multiplying both sides of the resulting inequality by $t>0$ and setting $Y(t):= t y(t)$, $g(t):=t G(t)$, we get
$$
Y^\prime(t)\le g(t)+y(t)+c_2y^{\alpha-1}(t)Y(t)\,,
$$
which entails
$$
Y(t)\le \beta+\int_0^th(s)Y(s){\rm d}s\,,\ \ t\in(0,\infty)\,,
$$
with
$$
\beta:=\int_0^\infty(g(s)+y(s)){\rm d}s\,,\ \ h(t):=c_2y^{\alpha-1}\,.
$$
Using Gronwall's lemma, we show
\be
Y(t)\le \beta \,{\rm exp}\left(\int_0^t h(s){\rm d}s\right)\,. 
\label{MoC}\ee
By assumption and \eqref{Co1} it follows that $\beta<\infty$ and $h\in L^1(0,\infty)$, so that the lemma follows from \eqref{MoC}.\par\hfill$\square$\par
{\em Proof of\, \theoref{nuovo}.} In \cite[Theorem 4.1]{GaSi}, it is shown the existence of a solution $(u,p,\xi,\omega,Q)$ to \eqref{eqpmodd}--\eqref{Q} in the class $\mathscr C_T$,\footnote{In \cite{GaSi} it is only shown $\nabla p\in L^2(0,\tau,L^2(\cald_R))$, for all $R\ge R_*$. However, proceeding exactly as in the proof of \cite[Lemma 4.3]{GaSi1}, one can demonstrate the stronger property $\nabla p\in L^2(0,\tau,L^2(\cald))$.} where $T$ is maximal, namely, either $T=\infty$, or else  there is $\{t_k\}\in (0,T)$ such that
\be\label{Si1}
\lim_{t_k\to T}\|D(u(t_k))\|_2=+\infty\,.
\ee
We shall show that, in fact, \eqref{Si1} cannot occur, provided the data satisfy \eqref{smalldata} for suitable $\delta>0$, thus implying that $(u,p,\xi,\omega,Q)$ exists for all times and is in $\mathscr C_\infty$. We begin to observe that, clearly,  $(u,p,\xi,\omega,Q)$ satisfies \eqref{casa} for all $t\in(0,T)$. Thus, using Cauchy-Schwarz on the right-hand side of \eqref{casa}$_1$ and integrating over $t\in (0,T)$ we get 
\be\ba{ll}\medskip\ds
\sup_{t\in[0,T]}\bigg( \varrho  \| u (t)\|_2^2 + m \left| \xi(t)
\right|^2 + \omega(t) \cdot I \cdot \omega(t)  \bigg) + \half\mu \int_0^T \| D(u(s)) \|_2^2\,{\rm d}s\\
\hspace*{3cm}\le \varrho  \| u_0\|_2^2 + m \left| \xi_0
\right|^2 + \omega_0 \cdot I \cdot \omega_0+c(\|\textup{F}\|_{L^2(0,\infty)}+\|\textup{M}\|_{L^2(0,\infty)})\,.
\label{Si2}\ea\ee
Moreover, again by Cauchy-Schwarz inequality and \lemmref{4.1},
\be\label{Si3}\ba{ll}\medskip
\big|m \omega \times \xi \cdot \dot{\xi}  +  \omega \times ({\sf I} \cdot \omega )\cdot \dot{\omega}+ {\sf F} \cdot \dot{\xi} + {\sf M} \cdot \dot{\omega}\big|\le \frac m2|\dot{\xi}|^2+\frac12\dot{\omega}\cdot{\sf I}\cdot\dot{\omega}+c\left(\|D(u)\|_2^4+|\textup{F}|^2+|\textup{M}|^2\right)\\
\ds \left|\int_{\cd} ( \omega \times \grad u_i\cdot\nabla u_i - \nabla(\omega \times u):\nabla u )\right|\le c\,\|D(u)\|_2^3\,.
\ea
\ee
Using the well-known trace inequality
$$
\|w\|_{2,\Sigma}\le c(\|w\|_2+\|w\|_2^\frac12\|\nabla w\|_2^\frac12)\,,\ \ w\in W^{1,2}(\cald)\,,
$$
along with \lemmref{4.1} and \lemmref{4.2}, we show
\be\label{Si4}\ba{rl}\medskip
\left|\Int{\Sigma}{}\left(n\cdot\nabla u\cdot\Phi-\half V\cdot n|\nabla u|^2\right)\right|&\le c\,(|\xi|+|\omega|)\|\nabla u\|^2_{2,\Sigma}\le c\, \big(\|D(u)\|_2^3+\|D(u)\|_2^2\|D^2u\|_2\big)\\
&\le c\, \big(\|D(u)\|_2^3+\|D(u)\|_2^4\big)+\frac14\|\Div T(u,p)\|_2^2\,.
\ea
\ee
Finally, employing the embedding inequality
$$
\|w\|_3\le c\,(\|\nabla w\|_2^\frac12\|w\|_2^\frac12+\|w\|_2)\,,\ \ w\in W^{1,2}(\cald)\,,
$$
H\"older inequality,  \lemmref{4.1} and \lemmref{4.2}, we show
\be\label{Si5}\ba{rl}\medskip\ds 
\left|\int_{\cd} u\cdot \grad u \cdot \dv T(u,p)\right|&\!\!\le c\,\|u\|_6\|D(u)\|_3\|\dv T(u,p)\|_2\le c\,\|D(u)\|_2(\|D^2u\|_2^\frac12\|D(u)\|_2^\frac12+\|D(u)\|_2)\|\Div T\|_2\\
&\!\!\le c\,(\|D(u)\|_2^4+\|D(u)\|_2^6)+\frac14\|\Div T\|_2^2\,.
\ea\ee
Using \eqref{Si3}--\eqref{Si5} in \eqref{casa}$_2$ we infer
\be\label{Si6}
\mu\frac{d}{dt}\| D(u)\|_{2}^2  +  \half m |\dot{\xi}|^2 + \half \dot{\omega} \cdot {\sf I} \cdot \dot{\omega} + \half \| \dv T(u,p)\|_{2}^2 
\le c\,(\| D(u)\|_{2}^3+\| D(u)\|_{2}^4+\| D(u)\|_{2}^6+|\textup{F}|^2+|\textup{M}|^2)\,,
\ee
and so, setting $y:=\|D(u)\|_2^2$, $G:=(c/\mu)(|\textup{F}|^2+|\textup{M}|^2)$ \eqref{Si6} furnishes, in particular, 
\be\label{Si7}
y^\prime(t)\le G(t)+ c\left(y^\frac32(t)+y^2(t)+y^3(t)\right)\,,\ \ t\in (0,T)\,,
\ee
Using multiple times Cauchy-Schwarz inequality, we show that \eqref{Si7} implies \eqref{1} with $\alpha=3$. Moreover, from \eqref{Si2}, we can find  $\delta>0$ such that if \eqref{smalldata} holds, then  assumption \eqref{2} of \lemmref{4.2} is satisfied.\footnote{Notice that from \eqref{casa}$_2$,  \eqref{Si3}--\eqref{Si5} and $u\in\mathscr C_T$ it follows that $d\|D(u)\|_2^2/dt\in L^1(0,\tau)$, all $\tau\in (0,T)$, that is, $\|D(u(t))\|_2^2$ is absolutely continuous in $t\in (0,T)$.} Thus, by that lemma, it follows that \be\sup_{t\in [0,T]}\|D(u(t))\|_2\le C\,,\label{COV}\ee which contradicts \eqref{Si1}. As a result, $T=\infty$ and therefore, by the second part of \lemmref{4.3}, we get
$$
\lim_{t\to\infty}\|D(u(t))\|_2=0\,.     
$$
The latter, in conjunction with \lemmref{4.1}, entails \eqref{asym}. We next observe that, from \eqref{Si2} and \lemmref{4.1}, it follows that
$$
\int_0^\infty\left(|{\xi}(t)|^2+|{\omega}(t)|^2+\|\nabla u(t)\|_2^2\right){\rm d}t\le C\,,
$$
while, by
integrating both sides of \eqref{Si6} over $t\in(0,\infty)$ and with the help of \lemmref{4.2} and \eqref{COV}, we get 
$$
\int_0^\infty\left(|\dot{\xi}(t)|^2+|\dot{\omega}(t)|^2+\|D^2u(t)\|_2^2\right){\rm d}t\le C\,,
$$
where the constant $C$ depends only on the data. Finally, from the latter, \eqref{Si5} and Sobolev inequality\footnote{Possibly, by modifying $p$ by adding a function of time.} we deduce
$$
\int_0^\infty\left(\|p(t)\|_6^2+\|\nabla p(t)\|_2^2\right){\rm d}t\le C\,,
$$
which completes the proof of the theorem.
\hfill$\square$
\Br  \theoref{nuovo} shows that, under the given assumptions on the data, the coupled system $\mathscr S$ will eventually go to a state of rest in the sense specified in \eqref{asym}. However, as also mentioned earlier on, we are not able to provide a rate of decay. In fact, we cannot apply the result stated in the last part of \lemmref{4.3} to the general case studied here, namely, a body of {\em arbitrary} shape. What prevents us from doing so is the presence of  $\|D(u)\|_2^3$ in \eqref{Si6} or, equivalently, $y^\frac32$ in \eqref{Si7}. However, if $\calb$ is a ball, that term does not occur. To show this, we notice that in this situation  the term $(\omega\times x\cdot\nabla u-\omega\times u)$ is no longer present in \eqref{eqpmodd}$_1$ as well as are not the terms $m\omega\times \xi$ and $\omega\times({\sf I}\cdot\omega)$ in \eqref{eqpmodd}$_{5}$, \eqref{eqpmodd}$_{6}$, respectively  (see, e.g., \cite{EMT}). Therefore, \eqref{casa}$_2$ becomes
$$  
\mu\frac{d}{dt}\| D(u)\|_{2}^2  +  m |\dot{\xi}|^2+\dot{\omega} \cdot {\sf I} \cdot \dot{\omega} + \| \dv T(u,p)\|_{2}^2 
=  -  {\sf F} \cdot \dot{\xi} - {\sf M} \cdot \dot{\omega} +    \varrho\ds\int_{\cd} (u-\xi)\cdot \grad u \cdot \dv T(u,p)\,.
$$
Arguing as in the proof of \theoref{nuovo}, one then shows
$$\ba{rl}\medskip
\left|  -  {\sf F} \cdot \dot{\xi} - {\sf M} \cdot \dot{\omega} +   \varrho\ds\int_{\cd} (u-\xi)\cdot \grad u \cdot \dv T(u,p)\right|
\le &\!\!\! \half \big( m |\dot{\xi}|^2+ \dot{\omega} \cdot {\sf I} \cdot \dot{\omega} + \|\dv T\|_2^2\big)\\&\!\!\!+\,c\,(|\textup{F}|^2+|\textup{M}|^2+\|D(u)\|_2^4+\|D(u)\|_2^6)\,. 
\ea 
$$
Consequently, combining the last two displayed relations and recalling that $\|D(u(t))\|_2$ is uniformly bounded in $t$ by the data, we deduce
$$
\frac{d}{dt}\| D(u)\|_{2}^2  
\le c\,(|\textup{F}|^2+|\textup{M}|^2+\| D(u)\|_{2}^4)\,,
$$
and by applying \lemmref{4.1} and \lemmref{4.3} we conclude
$$
\|u(t)\|_6+\|\nabla u(t)\|_2+|\xi(t)|+|\omega(t)|=O(t^{-\frac12})\ \ \mbox{as $t\to\infty$}\,,
$$
\label{rem:4.1}
\end{rm}\end{remark}
provided $t^\frac12{\rm F},\, t^\frac12{\rm M}\in L^2(0,\infty)$.
\medskip\par
{\bf Acknowledgment}. I would like to thank Professor Toshiaki Hishida, for bringing the problem to my attention. I also thank Mr. Jan A. Wein for several helpful conversations.

\end{document}